\documentclass[11pt,a4paper,reqno]{amsart}
\usepackage[margin=1in]{geometry}
\usepackage[T1]{fontenc}
\usepackage[utf8]{inputenc}
\usepackage{lmodern}
\usepackage{amssymb}
\usepackage{microtype}
\usepackage{enumitem}
\usepackage{hyperref}

\hypersetup{colorlinks=true,linkcolor=blue,citecolor=blue,urlcolor=blue,pdftitle={Box dimension of spiral trajectories near a degenerate focus of type (m,n)}}

\newtheorem{theorem}{Theorem}[section]
\newtheorem{lemma}[theorem]{Lemma}
\newtheorem{proposition}[theorem]{Proposition}
\newtheorem{corollary}[theorem]{Corollary}
\theoremstyle{definition}
\newtheorem{definition}[theorem]{Definition}
\newtheorem{remark}[theorem]{Remark}

\newcommand{\Cs}{\operatorname{Cs}}
\newcommand{\Sn}{\operatorname{Sn}}
\newcommand{\R}{\mathbb R}
\newcommand{\cT}{\mathcal T}
\newcommand{\cR}{\mathcal R}
\newcommand{\ldimB}{\underline{\dim}_{\mathrm B}}
\newcommand{\udimB}{\overline{\dim}_{\mathrm B}}
\newcommand{\dimB}{\dim_{\mathrm B}}
\newcommand{\mean}[1]{\left\langle #1\right\rangle}

\makeatletter
\newcommand{\affiliation}[1]{\gdef\@affiliation{#1}}
\let\@affiliation\@empty
\g@addto@macro\@setauthors{%
  \ifx\@affiliation\@empty\else
    \baselineskip13\p@
    \vtop{\centering{\footnotesize\itshape\@affiliation\@@par}%
      \global\dimen@i\prevdepth}\prevdepth\dimen@i
  \fi}
\makeatother

\begin{document}

\title{Box dimension of spiral trajectories near a degenerate focus of type $(m,n)$}
\author{Domagoj Vlah}
\affiliation{University of Zagreb Faculty of Electrical Engineering and Computing}

\begin{abstract}
We consider the family of planar polynomial vector fields with a degenerate focus of type $(m,n)$ at the origin studied by Huzak, Vlah, \v Zubrini\'c and \v Zupanovi\'c, and prove their conjectured formula for the box dimension of its spiral trajectories for all odd $m\ge n\ge1$ and all integers $k\ge1$. In generalized trigonometric (weighted polar) coordinates the radial equation can be integrated explicitly, and the two coordinates of a trajectory decay like $\varphi^{-p}$ and $\varphi^{-q}$, respectively, where
\[
p=\frac1{2mk},\qquad q=\frac1{2nk}.
\]
Decomposing a trajectory $\Gamma$ into successive turns, we prove the two-sided estimate
\[
N_\varepsilon(\Gamma)\asymp
\varepsilon^{-\frac{2+q-p}{1+q}}
\]
for the minimal number of sets of diameter at most $\varepsilon$ covering $\Gamma$. The upper bound combines the length $O(j^{-p})$ of the $j$th turn with a cover of the inner part of the spiral by a rectangle. For the lower bound we use a fixed sector in the first quadrant in which the turns are graphs over the horizontal axis; there, consecutive graphs are at vertical distance of order $j^{-q-1}$. Consequently $\Gamma$ is Minkowski nondegenerate and
\[
\dim_{\mathrm B}\Gamma
=\frac{2+q-p}{1+q}
=2-\frac{1+\frac nm}{1+2nk}.
\]
The argument adapts the covering and separation estimates of Burrell, Falconer and Fraser for elliptical polynomial spirals, but does not require a bi-Lipschitz comparison with such a spiral; for $m=n$ it gives a new proof of the known isotropic case. In addition, for a scalar periodic equation satisfied by $z=r^{2mn}$ along a trajectory, we construct a finite-order normal form with constant coefficients $g_j$ and show that if $g_1=\cdots=g_{\ell-1}=0$ and $g_\ell<0$, then $r(\varphi)\asymp\varphi^{-1/(2mn\ell)}$.
\end{abstract}

\maketitle

\section{Introduction}

Throughout, $\ldimB$ and $\udimB$ denote the lower and upper box dimension, and $\dimB$ denotes their common value when it exists. For a bounded set $E\subset\mathbb R^2$, $N_\varepsilon(E)$ is the least number of sets of diameter at most $\varepsilon$ needed to cover $E$. For positive functions $f,g$ we write $f\asymp g$ if $c\le f/g\le C$ for some constants $0<c\le C$.

Box dimensions of spirals were computed in \cite{DMT1983,Tricot1995}. \v Zubrini\'c and \v Zupanovi\'c \cite{ZZ2005,ZZ2008} initiated the fractal analysis of spiral trajectories of planar vector fields and related the box dimension of a spiral trajectory to the asymptotics of the Poincar\'e map. Fractal analysis of nilpotent and degenerate singularities and of slow-fast systems can be found in \cite{HDHVZ2021,HCV2021,HVZZ2023}, and spiral trajectories near pseudo foci of piecewise analytic systems are treated in \cite{Crnkovic2025}; see also the survey \cite{HRZ2025}.

The family \eqref{eq:system} of planar polynomial vector fields, where $m,n$ are odd and $k\ge1$, first appeared in the preprint \cite{VZZ2017}\footnote{In \cite[Theorem~3 and Remark~4]{VZZ2017}, the corresponding formula for $m>n$ contains an error: its denominator is $1+2km$ instead of $1+2kn$. The correct expression is the one conjectured in \cite{HVZZ2023} and proved in Theorem~\ref{thm:main} below.} and was studied by Huzak, Vlah, \v Zubrini\'c and \v Zupanovi\'c in \cite{HVZZ2023}. The origin is a focus, and for $(m,n)\neq(1,1)$ the linear part of the field at the origin is nilpotent or zero. For $m=n$, every spiral trajectory $\Gamma$ has box dimension $2-2/(1+2nk)$ and is Minkowski nondegenerate; this was announced in \cite[Theorem~1]{VZZ2017} and proved in \cite[Theorem~1]{HVZZ2023}. For $m\ge n$, the authors of \cite{HVZZ2023} conjectured, on the basis of numerical experiments, that every spiral trajectory is Minkowski nondegenerate and
\[
\dimB\Gamma=2-\frac{1+\frac nm}{1+2nk};
\]
for $m>n$, the dimension formula is also stated as a conjecture in the survey \cite[Section~2.2]{HRZ2025}. The main result of this paper, Theorem~\ref{thm:main} together with Corollary~\ref{cor:unstable}, proves this conjecture for all odd $m\ge n\ge1$ and $k\ge1$. For $m=n$ it gives a new proof of the box-dimension and nondegeneracy statements of \cite[Theorem~1]{HVZZ2023} for $k\ge1$; see Remark~\ref{rem:isotropic}.

The conjectured value is the box dimension of the elliptical polynomial spiral $S_{p,q}$ with $p=1/(2mk)$ and $q=1/(2nk)$, computed by Burrell, Falconer and Fraser \cite{BFF2022}; this was the heuristic behind the conjecture in \cite{HVZZ2023}. To deduce the formula from \cite{BFF2022}, one would need a bi-Lipschitz map between a trajectory of \eqref{eq:system} and $S_{p,q}$. Such a comparison is not automatic near the coordinate axes, where the oval $u^{2m}+v^{2n}=1$ traced by $(\Cs,\Sn)$ can be flat (for $(m,n)\neq(1,1)$ its curvature vanishes at $(0,\pm1)$), and we do not use it. Instead, the proof of Theorem~\ref{thm:main} works directly with the turns of the trajectory and adapts the covering and separation argument used for $S_{p,q}$ in \cite[Section~4]{BFF2022}: the first turns are covered using their lengths, the remaining ones by a rectangle, and the lower bound comes from the separation of consecutive turns. For the trajectories of \eqref{eq:system} the required estimates follow from the relation \eqref{eq:turns} between consecutive turns. The $j$th turn has length $O(j^{-p})$, and in a fixed sector of the first quadrant consecutive turns are graphs at vertical distance of order $j^{-q-1}$, which we show by implicit differentiation. Comparing this distance with $\varepsilon$ gives the critical turn number $j\asymp\varepsilon^{-1/(1+q)}$ and the exponent $(2+q-p)/(1+q)$.

Weighted polar coordinates, expansions of the Poincar\'e map, periodic near-identity transformations and generalized focal (Lyapunov) values are classical tools in the study of monodromic and nilpotent singularities; see, e.g., \cite{GarciaMaza2010,HanRomanovski2012,Sanders2007}. In Section~\ref{sec:nf}, which is independent of Theorem~\ref{thm:main}, we use them for general trajectories. If the radial equation along a trajectory has a finite expansion in powers of $r^{2mn}$, then $z=r^{2mn}$ satisfies a scalar $T$-periodic equation $dz/d\varphi=\sum_{j=1}^N a_j(\varphi)z^{j+1}+O(z^{N+2})$, and we construct a finite-order normal form for it with constant coefficients $g_j$ (Theorem~\ref{thm:nf}). In general the $g_j$ are not the averages of the $a_j$ (Remark~\ref{rem:averages}). If the first nonzero coefficient is $g_\ell<0$, then
\[
r(\varphi)\asymp\varphi^{-1/(2mn\ell)}
\]
(Theorem~\ref{thm:decay}). To obtain the box dimension of a general trajectory from this decay rate, our argument needs additional geometric information; in Corollary~\ref{cor:dim} this is a bi-Lipschitz comparison with $S_{p,q}$. For the model \eqref{eq:system}, the explicit turn structure used in Theorem~\ref{thm:main} makes such a hypothesis unnecessary.

The paper is organized as follows. Section~\ref{sec:prelim} recalls the generalized trigonometric functions and weighted polar coordinates. In Section~\ref{sec:model} we reduce \eqref{eq:system} to weighted polar coordinates and prove Theorem~\ref{thm:main}; Section~\ref{subsec:transfer} contains a transfer lemma for spirals that are piecewise bi-Lipschitz equivalent to $S_{p,q}$. Section~\ref{sec:nf} contains the normal form and the decay estimate, and Section~\ref{sec:concl} some concluding remarks.

\section{Generalized trigonometric functions and weighted polar coordinates}\label{sec:prelim}

Let $m,n\in\mathbb N$ be odd. The generalized trigonometric functions $\Cs$ and $\Sn$ are the solution of
\begin{equation}\label{eq:cssn}
\Cs'(\varphi)=-n\Sn(\varphi)^{2n-1},
\qquad
\Sn'(\varphi)=m\Cs(\varphi)^{2m-1},
\end{equation}
with initial conditions $\Cs(0)=1$, $\Sn(0)=0$. Differentiating shows that
\begin{equation}\label{eq:oval}
\Cs(\varphi)^{2m}+\Sn(\varphi)^{2n}=1.
\end{equation}
Hence the solution stays on the oval
\[
\{(u,v)\in\mathbb R^2:u^{2m}+v^{2n}=1\},
\]
a smooth closed curve on which the vector field $(u,v)\mapsto(-nv^{2n-1},mu^{2m-1})$ of \eqref{eq:cssn} does not vanish. Consequently the solution is defined and $C^\infty$ on $\R$, it is periodic, and it traverses the oval once in each period; we denote by $T>0$ its minimal period. For $m=1$ or $n=1$, functions of this kind were introduced, in a different normalization, by Lyapunov \cite{Lyapunov1893}. These functions and the associated weighted polar coordinates are standard tools in the study of degenerate and nilpotent singularities; see, e.g., \cite{HanRomanovski2012,HVZZ2023}.

We use the weighted polar coordinates (called $(n,m)$-polar coordinates in \cite{HVZZ2023})
\begin{equation}\label{eq:polar}
(x,y)=\bigl(r^n\Cs(\varphi),\ r^m\Sn(\varphi)\bigr).
\end{equation}
For $m=n=1$ these are the usual polar coordinates. By \eqref{eq:oval}, $x^{2m}+y^{2n}=r^{2mn}$, and a direct computation using \eqref{eq:cssn} shows that the Jacobian determinant of \eqref{eq:polar} with respect to $(r,\varphi)$ equals $mn\,r^{m+n-1}$. Thus every point of $\R^2\setminus\{0\}$ has coordinates $(r,\varphi)$ with $r>0$ uniquely determined and $\varphi$ uniquely determined modulo $T$, and along a continuous curve in $\R^2\setminus\{0\}$ the angle $\varphi$ can be chosen continuously.

\begin{definition}[Small spiral trajectory]\label{def:spiral}
A \emph{small spiral trajectory} of a planar vector field is a tail
\[
\Gamma=\{(x(\varphi),y(\varphi)):\varphi\ge \varphi_0\}
\]
of one of its trajectories, contained in $\R^2\setminus\{0\}$ and written in the coordinates \eqref{eq:polar}, such that $\varphi$ is strictly increasing along the trajectory, the map $\varphi\mapsto (x(\varphi),y(\varphi))$ is $C^1$, and $r(\varphi)\to 0$ as $\varphi\to\infty$. In particular, the image of every compact parameter interval is a compact $C^1$ arc.
\end{definition}

\begin{lemma}[The primitive of $\omega$]\label{lem:omega}
Let $m,n\in\mathbb N$ be odd, and let $\Cs,\Sn$ be defined by \eqref{eq:cssn}. Then the function
\[
\omega(\varphi):=\Sn(\varphi)^{n-1}\Cs(\varphi)^{m-1}
\]
is nonnegative and $T$-periodic, and
\[
\int_0^T \omega(\tau)\,d\tau>0.
\]
Consequently there exist a bounded $T$-periodic function $Q$ and a constant $\kappa>0$ such that
\[
\int_0^{\varphi}\omega(\tau)\,d\tau = \kappa \varphi + Q(\varphi)
\qquad (\varphi\in\R).
\]
\end{lemma}

\begin{proof}
Since $m$ and $n$ are odd, the exponents $m-1$ and $n-1$ are even, hence $\omega\ge 0$, and $\omega$ is $T$-periodic because $\Cs$ and $\Sn$ are. For small $\varphi>0$ we have $\Sn(\varphi)>0$ and $\Cs(\varphi)$ close to $1$, so $\omega>0$ on some interval; since $\omega$ is continuous and nonnegative, $\int_0^T\omega>0$. Set
\[
\kappa:=\frac1T\int_0^T \omega(\tau)\,d\tau>0,
\qquad
Q(\varphi):=\int_0^{\varphi}\omega(\tau)\,d\tau-\kappa\varphi .
\]
By the periodicity of $\omega$,
\[
Q(\varphi+T)-Q(\varphi)
=
\int_{\varphi}^{\varphi+T}\omega(\tau)\,d\tau-\kappa T
=
\int_0^T\omega(\tau)\,d\tau-\kappa T
=0,
\]
so $Q$ is $T$-periodic and, being continuous, bounded. In fact $\int_0^T\omega(\tau)\,d\tau=2\pi/(mn)$ \cite[Lemma~2]{HVZZ2023}, but we only need its positivity.
\end{proof}

\section{The model family}\label{sec:model}

\subsection{Reduction to weighted polar coordinates}

The following family was introduced in \cite{HVZZ2023}. We recall the reduction of its stable branch to weighted polar coordinates (cf.\ \cite[eq.~(42) and Theorem~2]{HVZZ2023}), because the explicit form of the trajectories is used in the proof of Theorem~\ref{thm:main}.

Fix odd integers $m\ge n\ge 1$ and an integer $k\ge 1$, and consider
\begin{equation}\label{eq:system}
\begin{aligned}
\dot x&=-n y^{2n-1}-n x^m y^{n-1}(x^{2m}+y^{2n})^k,\\
\dot y&=\phantom{-}m x^{2m-1}-m x^{m-1} y^n (x^{2m}+y^{2n})^k.
\end{aligned}
\end{equation}
The minus signs in front of the terms containing $(x^{2m}+y^{2n})^k$ make the origin a stable focus; see \eqref{eq:radial}. The field with the opposite signs is treated in Corollary~\ref{cor:unstable}.

\begin{proposition}[Weighted polar form of \eqref{eq:system}]
Along every nonconstant trajectory of \eqref{eq:system}, written in the coordinates \eqref{eq:polar}, one has
\begin{equation}\label{eq:phidot}
\dot\varphi=r^{2mn-m-n},
\end{equation}
\begin{equation}\label{eq:rdot}
\dot r=-\Sn(\varphi)^{n-1}\Cs(\varphi)^{m-1}r^{2mn-m-n+1+2mnk}.
\end{equation}
In particular, $\varphi$ is strictly increasing, and
\begin{equation}\label{eq:radial}
\frac{dr}{d\varphi}=-\Sn(\varphi)^{n-1}\Cs(\varphi)^{m-1}r^{2mnk+1}.
\end{equation}
Along every small spiral trajectory of \eqref{eq:system} there is a constant $C_0\in\R$ such that
\begin{equation}\label{eq:rint}
r(\varphi)=\Bigl(2mnk\int_0^{\varphi}\Sn(\tau)^{n-1}\Cs(\tau)^{m-1}\,d\tau + C_0\Bigr)^{-1/(2mnk)}.
\end{equation}
Furthermore, since $m$ and $n$ are odd, there exist $K>0$ and a bounded $T$-periodic function $\widetilde Q$ such that
\begin{equation}\label{eq:rnorm}
r(\varphi)=K\bigl(\varphi+\widetilde Q(\varphi)\bigr)^{-1/(2mnk)}.
\end{equation}
Hence the trajectory has the explicit parametrization
\begin{equation}\label{eq:traj}
\begin{aligned}
x(\varphi)&=K^n\Cs(\varphi)\bigl(\varphi+\widetilde Q(\varphi)\bigr)^{-1/(2mk)},\\
y(\varphi)&=K^m\Sn(\varphi)\bigl(\varphi+\widetilde Q(\varphi)\bigr)^{-1/(2nk)}.
\end{aligned}
\end{equation}
\end{proposition}

\begin{proof}
Write
\[
C:=\Cs(\varphi),\qquad S:=\Sn(\varphi).
\]
Differentiating \eqref{eq:polar} and using \eqref{eq:cssn} yields
\begin{equation}\label{eq:xydot}
\begin{aligned}
\dot x&=n r^{n-1}\dot r\,C-n r^n S^{2n-1}\dot\varphi,\\
\dot y&=m r^{m-1}\dot r\,S+m r^m C^{2m-1}\dot\varphi.
\end{aligned}
\end{equation}
On the other hand, \eqref{eq:oval} implies
\[
x^{2m}+y^{2n}=r^{2mn}(C^{2m}+S^{2n})=r^{2mn},
\]
so \eqref{eq:system} becomes
\begin{equation}\label{eq:system-polar}
\begin{aligned}
\dot x&=-n r^{2mn-m}S^{2n-1}-n r^{2mn-m+2mnk}C^mS^{n-1},\\
\dot y&=\phantom{-}m r^{2mn-n}C^{2m-1}-m r^{2mn-n+2mnk}C^{m-1}S^n.
\end{aligned}
\end{equation}
Dividing the first identity in \eqref{eq:xydot} by $n r^{n-1}$ and the second by $m r^{m-1}$, and comparing with \eqref{eq:system-polar}, we obtain the linear system
\begin{equation}
\begin{aligned}
C\,\dot r-rS^{2n-1}\dot\varphi
&=
-r^{2mn-m-n+1}\bigl(S^{2n-1}+r^{2mnk}C^mS^{n-1}\bigr),\\
S\,\dot r+rC^{2m-1}\dot\varphi
&=
\phantom{-}r^{2mn-m-n+1}\bigl(C^{2m-1}-r^{2mnk}C^{m-1}S^n\bigr)
\end{aligned}
\end{equation}
for $(\dot r,\dot\varphi)$. Its determinant is
\[
r\bigl(C^{2m}+S^{2n}\bigr)=r.
\]
By Cramer's rule,
\begin{align*}
\dot\varphi
&=
r^{2mn-m-n}
\frac{C\bigl(C^{2m-1}-r^{2mnk}C^{m-1}S^n\bigr)+S\bigl(S^{2n-1}+r^{2mnk}C^mS^{n-1}\bigr)}
{C^{2m}+S^{2n}}\\
&=
r^{2mn-m-n},
\end{align*}
because the mixed terms cancel and $C^{2m}+S^{2n}=1$. Likewise,
\begin{align*}
\dot r
&=
r^{2mn-m-n+1}
\frac{-\bigl(S^{2n-1}+r^{2mnk}C^mS^{n-1}\bigr)C^{2m-1}
+S^{2n-1}\bigl(C^{2m-1}-r^{2mnk}C^{m-1}S^n\bigr)}
{C^{2m}+S^{2n}}\\
&=
-\Sn(\varphi)^{n-1}\Cs(\varphi)^{m-1}r^{2mn-m-n+1+2mnk}.
\end{align*}
Therefore $\dot\varphi>0$, and $\varphi$ is strictly increasing. Dividing \eqref{eq:rdot} by \eqref{eq:phidot} yields \eqref{eq:radial}. Equivalently,
\[
\frac{d}{d\varphi}\bigl(r^{-2mnk}\bigr)=2mnk\,\Sn(\varphi)^{n-1}\Cs(\varphi)^{m-1},
\]
and integrating gives \eqref{eq:rint} for some constant $C_0$. By Lemma~\ref{lem:omega}, which uses that $m$ and $n$ are odd,
\[
2mnk\int_0^{\varphi}\Sn(\tau)^{n-1}\Cs(\tau)^{m-1}\,d\tau + C_0
=2mnk\kappa\bigl(\varphi+\widetilde Q(\varphi)\bigr),
\qquad
\widetilde Q:=\frac{Q}{\kappa}+\frac{C_0}{2mnk\kappa},
\]
which is \eqref{eq:rnorm} with $K=(2mnk\kappa)^{-1/(2mnk)}$. Substituting \eqref{eq:rnorm} into \eqref{eq:polar} gives \eqref{eq:traj}.
\end{proof}

\begin{proposition}[Asymptotics of the coordinates]
For every small spiral trajectory of \eqref{eq:system}, one has
\[
x(\varphi)=A(\varphi)\bigl(\varphi+\widetilde Q(\varphi)\bigr)^{-1/(2mk)},
\qquad
y(\varphi)=B(\varphi)\bigl(\varphi+\widetilde Q(\varphi)\bigr)^{-1/(2nk)},
\]
where $\widetilde Q$ is bounded and periodic, and
\[
A(\varphi)=K^n\Cs(\varphi),\qquad B(\varphi)=K^m\Sn(\varphi)
\]
are bounded $T$-periodic functions that do not vanish identically. In particular,
\[
x(\varphi)=A(\varphi)\varphi^{-1/(2mk)}+o\bigl(\varphi^{-1/(2mk)}\bigr),
\qquad
y(\varphi)=B(\varphi)\varphi^{-1/(2nk)}+o\bigl(\varphi^{-1/(2nk)}\bigr).
\]
\end{proposition}

\begin{proof}
The first two formulas are \eqref{eq:traj}. Since $\widetilde Q$ is bounded,
\[
\bigl(\varphi+\widetilde Q(\varphi)\bigr)^{-\alpha}
=
\varphi^{-\alpha}\Bigl(1+\frac{\widetilde Q(\varphi)}{\varphi}\Bigr)^{-\alpha}
=
\varphi^{-\alpha}+O\bigl(\varphi^{-\alpha-1}\bigr)
\qquad (\varphi\to\infty)
\]
for every fixed $\alpha>0$. Applying this with $\alpha=1/(2mk)$ and $\alpha=1/(2nk)$ gives the asymptotic formulas for $x(\varphi)$ and $y(\varphi)$.
\end{proof}

\subsection{Box dimension of the spiral trajectories}

\begin{theorem}[Box dimension of the model spiral]\label{thm:main}
Let $m\ge n\ge1$ be odd and $k\ge1$. For every small spiral trajectory $\Gamma$ of the stable model field \eqref{eq:system}, set
\[
p:=\frac1{2mk},\qquad q:=\frac1{2nk},\qquad
d:=\frac{2+q-p}{1+q}.
\]
There exist constants $c,C,\varepsilon_0>0$, depending on $\Gamma$, such that
\begin{equation}\label{eq:cover}
c\varepsilon^{-d}\le N_\varepsilon(\Gamma)\le C\varepsilon^{-d}
\qquad (0<\varepsilon<\varepsilon_0).
\end{equation}
Consequently,
\begin{equation}\label{eq:dim}
\dimB\Gamma=d
=2-\frac{1+\frac nm}{1+2nk}.
\end{equation}
Moreover, if $\Gamma_\varepsilon$ denotes the Euclidean $\varepsilon$-neighbourhood of $\Gamma$ and $|\cdot|$ denotes planar Lebesgue measure, then
\begin{equation}\label{eq:mink}
0<\liminf_{\varepsilon\to0^+}
\frac{|\Gamma_\varepsilon|}{\varepsilon^{2-d}}
\le
\limsup_{\varepsilon\to0^+}
\frac{|\Gamma_\varepsilon|}{\varepsilon^{2-d}}
<\infty,
\end{equation}
that is, $\Gamma$ is Minkowski nondegenerate in dimension $d$.
\end{theorem}

\begin{proof}
Throughout the proof, $c$ and $C$ denote positive constants that do not depend on $\varepsilon$ or on the turn index; their values may change from line to line. Put
\[
\omega(\varphi)=\Sn(\varphi)^{n-1}\Cs(\varphi)^{m-1},
\qquad
H(\varphi)=2mnk\int_0^\varphi\omega(\tau)\,d\tau+C_0.
\]
By \eqref{eq:rint}, $H=r^{-2mnk}>0$ along $\Gamma$, and \eqref{eq:polar} gives
\begin{equation}
x(\varphi)=\Cs(\varphi)H(\varphi)^{-p},
\qquad
y(\varphi)=\Sn(\varphi)H(\varphi)^{-q}.
\end{equation}
Let
\[
D:=2mnk\int_0^T\omega(\tau)\,d\tau>0,
\qquad h(s):=H(s),\quad 0\le s\le T.
\]
By the periodicity of $\omega$, consecutive turns are related by
\begin{equation}\label{eq:turns}
H(jT+s)=jD+h(s)
\qquad (j\in\mathbb N,\ 0\le s\le T).
\end{equation}
Since $h$ is bounded on $[0,T]$, there are $j_0\in\mathbb N$ with $j_0T\ge\varphi_0$ and constants $0<c_H\le C_H$ such that
\begin{equation}\label{eq:Hu}
c_Hu\le uD+h(s)\le C_Hu
\qquad (u\ge j_0,\ s\in[0,T]).
\end{equation}
The arc $\{(x(\varphi),y(\varphi)):\varphi_0\le\varphi\le j_0T\}$ is a compact $C^1$ arc, so its covering number is $O(\varepsilon^{-1})$; since $d>1$, it does not affect \eqref{eq:cover}, and we disregard it. For integers $j\ge j_0$ the $j$th turn is
\[
\Gamma_j:=\{(X_j(s),Y_j(s)):0\le s\le T\},
\]
where
\[
X_j(s)=\Cs(s)(jD+h(s))^{-p},
\qquad
Y_j(s)=\Sn(s)(jD+h(s))^{-q}.
\]

\emph{Upper bound.} The functions $\Cs,\Sn$ and their derivatives are bounded, and so is $h'=2mnk\omega$. Differentiating $X_j,Y_j$ and using \eqref{eq:Hu} and $q\ge p$, we get
\[
|X_j'(s)|+|Y_j'(s)|\le Cj^{-p}
\qquad (j\ge j_0,\ 0\le s\le T).
\]
Hence $\Gamma_j$ has length at most $Cj^{-p}$ and can be covered by at most $1+Cj^{-p}/\varepsilon$ sets of diameter at most $\varepsilon$. Let
\[
J:=\left\lceil\varepsilon^{-1/(1+q)}\right\rceil.
\]
Since $0<p<1$, we have $J^p\le\varepsilon^{-1}$ for small $\varepsilon$, and
\[
d=1+\frac{1-p}{1+q}.
\]
Therefore the turns with $j_0\le j\le J$ can be covered by at most
\[
\sum_{j=j_0}^{J}\Bigl(1+\frac{Cj^{-p}}{\varepsilon}\Bigr)
\le J+C\varepsilon^{-1}J^{1-p}
\le C\varepsilon^{-1}J^{1-p}
\le C\varepsilon^{-d}
\]
sets of diameter at most $\varepsilon$. By \eqref{eq:Hu}, the remaining part $\{(x(\varphi),y(\varphi)):\varphi\ge (J+1)T\}$ of $\Gamma$ lies in a rectangle with side lengths at most $CJ^{-p}$ and $CJ^{-q}$. Covering this rectangle by squares of side $\varepsilon/\sqrt2$ requires at most
\[
C\left(1+\frac{J^{-p}}\varepsilon\right)
 \left(1+\frac{J^{-q}}\varepsilon\right)
\le C\varepsilon^{-2}J^{-p-q}
\le C\varepsilon^{-d}
\]
squares, because $J^{-p}\ge J^{-q}\ge\varepsilon$ for small $\varepsilon$ and $2-(p+q)/(1+q)=d$. This proves $N_\varepsilon(\Gamma)\le C\varepsilon^{-d}$.

\emph{Lower bound.} Choose compact intervals $I=[a,b]$ and $I'=[a',b']$ with $a'<a<b<b'$, both contained in an open subinterval of $(0,T)$ on which $\Cs>0$ and $\Sn>0$, that is, in the parameter interval of the arc of the oval lying in the open first quadrant. By compactness, there is $c_1>0$ such that on $I'$
\[
\Cs,\Sn\ge c_1,
\qquad
\Cs'\le-c_1,
\qquad
\Sn'\ge c_1.
\]
For real $u\ge j_0$ set
\[
H_u(s)=uD+h(s),\qquad
X_u(s)=\Cs(s)H_u(s)^{-p},\qquad
Y_u(s)=\Sn(s)H_u(s)^{-q}
\quad(s\in I').
\]
Since $h'\ge0$,
\[
\partial_sX_u
=H_u^{-p}\left(\Cs'-\frac{p\Cs h'}{H_u}\right)<0
\quad\text{on }I',
\]
and, for sufficiently large $u$,
\[
\partial_sY_u
=H_u^{-q}\left(\Sn'-\frac{q\Sn h'}{H_u}\right)>0
\quad\text{on }I'.
\]
Thus, for large $u$, the arc $\{(X_u(s),Y_u(s)):s\in I\}$ is the graph $y=G_u(x)$ of a function on the interval
\[
I_u=[X_u(b),X_u(a)],
\qquad |I_u|\asymp u^{-p}.
\]
Both endpoints of $I_u$ decrease with $u$, because
\[
\partial_uX_u=-pD\Cs H_u^{-p-1}<0.
\]
For $x\in I_u$ let $s(u,x)\in I$ be the solution of $X_u(s)=x$. Since $\partial_sX_u<0$ on $I'$ and $s(u,x)$ lies in the interior of $I'$, the implicit function theorem shows that $s(u,x)$ and $G_u(x)=Y_u(s(u,x))$ are $C^1$ functions of $u$ on every interval of values of $u$ for which $x\in I_u$, and
\[
\partial_uG_u(x)
=D H_u^{-q-1}
\left[
-q\Sn
+p\Cs\,
\frac{\Sn'-q\Sn h'/H_u}{\Cs'-p\Cs h'/H_u}
\right],
\]
where all functions of $s$ are evaluated at $s=s(u,x)$. The second term in brackets is negative, while $\Sn\ge c_1$. The denominator is bounded away from zero on $I'$, and all remaining factors are bounded. Together with \eqref{eq:Hu}, this gives constants $c_2,C_3>0$ and $u_0\ge j_0$ such that
\begin{equation}\label{eq:dG}
c_2u^{-q-1}
\le-\partial_uG_u(x)
\le C_3u^{-q-1}
\end{equation}
for every $u\ge u_0$ and $x\in I_u$. For integers $j\ge u_0$, the graph of $G_j$ is a subarc of $\Gamma_j\subset\Gamma$.

Let $u_0\le j<j'$ be integers and $x\in I_j\cap I_{j'}$. Since both endpoints of $I_u$ decrease with $u$, we have $x\in I_u$ for every $u\in[j,j']$, and integrating \eqref{eq:dG} over $[j,j']$ gives
\begin{equation}\label{eq:gap}
G_j(x)-G_{j'}(x)
\ge c_2(j'-j)(j')^{-q-1}.
\end{equation}
Choose $A>0$ so small that $c_2A^{-(q+1)}\ge8$, and put
\[
J_\varepsilon=\left\lfloor A\varepsilon^{-1/(1+q)}\right\rfloor.
\]
For sufficiently small $\varepsilon$ we have $J_\varepsilon\ge 2u_0$, and \eqref{eq:gap} shows that if $u_0\le j<j'\le J_\varepsilon$ and $x\in I_j\cap I_{j'}$, then
\[
G_j(x)-G_{j'}(x)\ge c_2J_\varepsilon^{-q-1}\ge c_2A^{-(q+1)}\varepsilon\ge 8\varepsilon.
\]

For every integer $j\in[u_0,J_\varepsilon]$, select the points of the graph of $G_j$ whose horizontal coordinate belongs to $4\varepsilon\mathbb Z\cap I_j$, and let $\mathcal P_\varepsilon$ be the set of all selected points. Since $|I_j|\ge cj^{-p}\ge cJ_\varepsilon^{-p}$ and $J_\varepsilon^{-p}/\varepsilon\to\infty$ as $\varepsilon\to0^+$, at least $c\varepsilon^{-1}j^{-p}$ points are selected on the $j$th graph for all sufficiently small $\varepsilon$. Therefore
\[
\#\mathcal P_\varepsilon
\ge c\varepsilon^{-1}
\sum_{u_0\le j\le J_\varepsilon}j^{-p}
\ge c\varepsilon^{-1}J_\varepsilon^{1-p}
\ge c\varepsilon^{-d}.
\]
Two distinct points of $\mathcal P_\varepsilon$ either have horizontal coordinates that differ by at least $4\varepsilon$, or have the same horizontal coordinate and lie on different graphs, in which case their vertical distance is at least $8\varepsilon$. Hence $\mathcal P_\varepsilon$ is $4\varepsilon$-separated, and a set of diameter at most $\varepsilon$ contains at most one of its points. Since $\mathcal P_\varepsilon\subset\Gamma$,
\[
N_\varepsilon(\Gamma)\ge\#\mathcal P_\varepsilon\ge c\varepsilon^{-d}.
\]
This proves \eqref{eq:cover} and therefore \eqref{eq:dim}.

\emph{Minkowski content.} The open disks of radius $\varepsilon$ centred at the points of $\mathcal P_\varepsilon$ are pairwise disjoint and contained in $\Gamma_\varepsilon$. Thus
\[
|\Gamma_\varepsilon|
\ge \pi\varepsilon^2\#\mathcal P_\varepsilon
\ge c\varepsilon^{2-d}.
\]
Conversely, if $U$ is a set of diameter at most $\varepsilon$ and $z\in U$, then the $\varepsilon$-neighbourhood of $U$ is contained in the disk of radius $2\varepsilon$ centred at $z$. Applying this to the sets of a minimal cover of $\Gamma$ by sets of diameter at most $\varepsilon$, and using the upper bound, we obtain
\[
|\Gamma_\varepsilon|\le 4\pi\varepsilon^2N_\varepsilon(\Gamma)
\le C\varepsilon^{2-d}.
\]
This is \eqref{eq:mink}.
\end{proof}

\begin{corollary}[The unstable case]\label{cor:unstable}
Let $F_+$ be the field obtained from \eqref{eq:system} by replacing the two minus signs in front of the terms containing $(x^{2m}+y^{2n})^k$ by plus signs. The conclusions of Theorem~\ref{thm:main} also hold for every trajectory tail of $F_+$ that tends to the origin in backward time.
\end{corollary}

\begin{proof}
Let $F_-$ denote the field \eqref{eq:system} and set $R(x,y)=(x,-y)$. Because $m$ and $n$ are odd, direct substitution gives
\[
-DR\,F_+(Rz)=F_-(z).
\]
Thus, if $z(t)$ solves $z'=F_+(z)$, then $w(t):=Rz(-t)$ solves $w'=F_-(w)$. Hence $R$ maps every trajectory tail of $F_+$ that tends to the origin in backward time onto a trajectory tail of $F_-$ that tends to the origin in forward time; by \eqref{eq:rint}, $\varphi\to\infty$ along such a tail, so it is a small spiral trajectory of $F_-$. Since $R$ is an isometry, it preserves covering numbers and the areas of $\varepsilon$-neighbourhoods, and the claim follows from Theorem~\ref{thm:main}.
\end{proof}

\begin{remark}[The case $m<n$]\label{rem:swap}
Denote the field \eqref{eq:system} by $F_{m,n}$ and let $S(x,y)=(y,-x)$. Since $m$ and $n$ are odd, a direct computation gives $DS\,F_{m,n}=F_{n,m}\circ S$, so the rotation $S$ maps trajectories of $F_{m,n}$ onto trajectories of $F_{n,m}$, preserving the direction of time. Hence, if $m<n$, Theorem~\ref{thm:main} applied to $F_{n,m}$ shows that every small spiral trajectory $\Gamma$ of $F_{m,n}$ is Minkowski nondegenerate and $\dimB\Gamma=2-(1+\frac mn)/(1+2mk)$. The assumption $m\ge n$ is therefore only a normalization.
\end{remark}

\subsection{A transfer lemma}\label{subsec:transfer}

Theorem~\ref{thm:main} does not use a comparison with an elliptical polynomial spiral. For trajectories without an explicit turn structure, however, the box dimension can be obtained from \cite{BFF2022} once a bi-Lipschitz comparison with such a spiral is available. This is the content of Lemma~\ref{lem:transfer}, which is used in Corollary~\ref{cor:dim}.

Write
\[
S_{p,q}:=\{(t^{-p}\cos t,\,t^{-q}\sin t): t\ge t_0\},
\qquad t_0>0,\quad 0<p\le q.
\]
Burrell, Falconer, and Fraser proved that
\begin{equation}\label{eq:bff}
\dimB S_{p,q}=\frac{2+q-p}{1+q}
\end{equation}
for $0<p<1$, while $\dimB S_{p,q}=1$ for $p\ge 1$; see \cite[Corollary~2.3]{BFF2022}, where the parameter range is $t>1$. For H\"older maps between such spirals, see \cite{BFF2022,CGV2026}. The lower cutoff $t_0>0$ is immaterial for box dimension: replacing one positive cutoff by another adds or removes only the image of a compact parameter interval, which is a rectifiable $C^1$ arc and hence has upper box dimension at most $1$. Since $(2+q-p)/(1+q)>1$ for $0<p<1$, this changes neither the upper box dimension nor, by the argument used for $\Sigma_1$ in the proof of Lemma~\ref{lem:transfer}, the lower box dimension.

\begin{definition}[Finite-tail bi-Lipschitz comparison]\label{def:comparison}
Let $0<p\le q$ with $p<1$, and let $\Gamma$ be a small spiral trajectory. We say that $\Gamma$ satisfies the \emph{finite-tail bi-Lipschitz comparison hypothesis with $S_{p,q}$} if there exist compact initial arcs $K_\Gamma\subset \Gamma$ and $K_S\subset S_{p,q}$ and decompositions
\[
\Gamma\setminus K_\Gamma=\Gamma_1\cup\cdots\cup\Gamma_M,
\qquad
S_{p,q}\setminus K_S=\Sigma_1\cup\cdots\cup\Sigma_M
\]
with the following properties.
\begin{enumerate}[label=\textnormal{(\alph*)},leftmargin=2em]
\item $K_\Gamma$ and $K_S$ are the images of compact initial intervals of the natural parameters $\varphi$ of $\Gamma$ and $t$ of $S_{p,q}$, respectively, and each $\Gamma_i$ and each $\Sigma_i$ is the image of an interval of the corresponding parameter; in particular, all pieces are connected.
\item For each $i\in\{1,\dots,M\}$ there is a bi-Lipschitz bijection $F_i:\Sigma_i\to\Gamma_i$.
\item $\Gamma_1$ and $\Sigma_1$ are the only pieces that are images of unbounded parameter intervals.
\item Distinct pieces of each decomposition meet at most at common boundary points.
\end{enumerate}
\end{definition}

\begin{lemma}[Transfer lemma]\label{lem:transfer}
Let $0<p\le q$ with $p<1$. If a small spiral trajectory $\Gamma$ satisfies the hypothesis of Definition~\ref{def:comparison} with $S_{p,q}$, then
\[
\ldimB \Gamma=\udimB \Gamma=\dimB \Gamma = \frac{2+q-p}{1+q}=2-\frac{p+q}{1+q}.
\]
\end{lemma}

\begin{proof}
Set $d:=(2+q-p)/(1+q)$. Then $d-1=(1-p)/(1+q)>0$, and $\ldimB S_{p,q}=\udimB S_{p,q}=d$ by \eqref{eq:bff}. We use two standard facts. First, if $E\subset\R^2$ is a compact rectifiable arc, then there is $L_E>0$ such that $N_\delta(E)\le L_E\delta^{-1}$ for all sufficiently small $\delta>0$; to see this, split $E$ into consecutive subarcs of length at most $\delta$. Second, if $F:E\to E'$ is a surjective Lipschitz map with Lipschitz constant $A>0$, then $N_\delta(E')\le N_{\delta/A}(E)$; in particular, lower and upper box dimensions are invariant under bi-Lipschitz bijections (see, e.g., \cite{Falconer2014}).

Let $\gamma(t):=(t^{-p}\cos t,\,t^{-q}\sin t)$. By (a) and (c), $K_S$ is the image under $\gamma$ of a compact interval, and for $i\ge2$ the piece $\Sigma_i$ is contained in $\gamma([a_i,b_i])$ for some compact interval $[a_i,b_i]\subset[t_0,\infty)$. Since $\gamma$ is smooth, these sets are compact rectifiable arcs. Hence the set
\[
C_S:=K_S\cup\Sigma_2\cup\cdots\cup\Sigma_M
\]
satisfies $N_\delta(C_S)\le L\delta^{-1}$ for all sufficiently small $\delta>0$ and some $L>0$. Fix $\eta\in(0,d-1)$. Since $\ldimB S_{p,q}=d$, there are $c_1,\delta_0>0$ such that $N_\delta(S_{p,q})\ge c_1\delta^{-(d-\eta)}$ for $0<\delta<\delta_0$. As $S_{p,q}=\Sigma_1\cup C_S$,
\[
N_\delta(\Sigma_1)\ge N_\delta(S_{p,q})-N_\delta(C_S)\ge c_1\delta^{-(d-\eta)}-L\delta^{-1}\ge \frac{c_1}{2}\,\delta^{-(d-\eta)}
\]
for all sufficiently small $\delta$, because $d-\eta>1$. Hence $\ldimB\Sigma_1\ge d-\eta$, and since $\eta$ was arbitrary, $\ldimB\Sigma_1\ge d$. On the other hand, $\udimB\Sigma_1\le\udimB S_{p,q}=d$. Therefore $\ldimB\Sigma_1=\udimB\Sigma_1=d$.

By (b), $F_1:\Sigma_1\to\Gamma_1$ is bi-Lipschitz, so $\ldimB\Gamma_1=\udimB\Gamma_1=d$. For $i\ge 2$, if $A_i$ is a Lipschitz constant of $F_i$, then
\[
N_\delta(\Gamma_i)\le N_{\delta/A_i}(\Sigma_i)\le L_iA_i\,\delta^{-1}
\]
for all sufficiently small $\delta>0$, where $L_i:=L_{\gamma([a_i,b_i])}$. By Definition~\ref{def:spiral}, $K_\Gamma$ is the image of a compact interval under the $C^1$ map $\varphi\mapsto (x(\varphi),y(\varphi))$, hence a compact rectifiable arc. Therefore the set
\[
C_\Gamma:=K_\Gamma\cup\Gamma_2\cup\cdots\cup\Gamma_M
\]
satisfies $N_\delta(C_\Gamma)\le L'\delta^{-1}$ for all sufficiently small $\delta$ and some $L'>0$. Since $\Gamma_1\subset \Gamma$, we have $\ldimB\Gamma\ge\ldimB\Gamma_1=d$. Finally, fix $\eta>0$. Because $\udimB\Gamma_1=d$, there are $A,\delta_1>0$ such that $N_\delta(\Gamma_1)\le A\delta^{-(d+\eta)}$ for $0<\delta<\delta_1$, and then
\[
N_\delta(\Gamma)
\le N_\delta(\Gamma_1)+N_\delta(C_\Gamma)
\le A\delta^{-(d+\eta)}+L'\delta^{-1}
\le 2A\delta^{-(d+\eta)}
\]
for all sufficiently small $\delta$, because $d+\eta>1$. Hence $\udimB\Gamma\le d+\eta$ for every $\eta>0$, and $\ldimB\Gamma=\udimB\Gamma=d$.
\end{proof}

\begin{remark}
The connectedness of the pieces and conditions (c) and (d) in Definition~\ref{def:comparison} are convenient but not essential. Suppose only that $\Gamma\setminus K_\Gamma=\Gamma_1\cup\cdots\cup\Gamma_M$ and $S_{p,q}\setminus K_S=\Sigma_1\cup\cdots\cup\Sigma_M$ for arbitrary sets $\Gamma_i,\Sigma_i$, with compact rectifiable arcs $K_\Gamma$ and $K_S$, and that there are bijections $F_i:\Sigma_i\to\Gamma_i$ such that $F_i$ and $F_i^{-1}$ are Lipschitz with constant at most $A$. Then $N_{A\delta}(\Sigma_i)\le N_\delta(\Gamma_i)\le N_{\delta/A}(\Sigma_i)$ for every $i$, and therefore
\[
\frac1M\bigl(N_{A\delta}(S_{p,q})-N_{A\delta}(K_S)\bigr)
\le N_\delta(\Gamma)
\le N_\delta(K_\Gamma)+M\,N_{\delta/A}(S_{p,q}).
\]
Since $N_\delta(K_\Gamma)$ and $N_\delta(K_S)$ are $O(\delta^{-1})$ and $(2+q-p)/(1+q)>1$, this again gives $\dimB\Gamma=(2+q-p)/(1+q)$. In particular, decompositions into finitely many families of sectors of the turns, which are not connected, are also admissible.
\end{remark}

\subsection{Remarks}

\begin{remark}[The isotropic case]\label{rem:isotropic}
For $m=n$ one has $p=q=1/(2nk)$, and \eqref{eq:dim} reads
\[
\dimB\Gamma=\frac{2}{1+q}=2-\frac{2}{1+2nk}.
\]
Together with Corollary~\ref{cor:unstable}, Theorem~\ref{thm:main} therefore gives a new proof of the box-dimension and Minkowski nondegeneracy statements of \cite[Theorem~1]{HVZZ2023} for $k\ge1$ (the comparability with a power spiral stated there follows directly from \eqref{eq:rnorm}). The proof in \cite{HVZZ2023} uses a bi-Lipschitz map between the trajectory and the power spiral $r=\varphi^{-1/(2nk)}$, whereas the proof of Theorem~\ref{thm:main} does not compare the trajectory with a model spiral.
\end{remark}

\begin{remark}[Minkowski measurability]\label{rem:measurability}
Theorem~\ref{thm:main} gives Minkowski nondegeneracy, that is, two-sided bounds for $|\Gamma_\varepsilon|/\varepsilon^{2-d}$. It does not show that this ratio converges as $\varepsilon\to0^+$, that is, that $\Gamma$ is Minkowski measurable; the covering argument determines $|\Gamma_\varepsilon|$ only up to multiplicative constants and cannot decide this question. For $m=n=1$, system \eqref{eq:system} reads $\dot r=-r^{2k+1}$, $\dot\varphi=1$ in polar coordinates, which is the normal form of a weak focus of order $k$. Up to the reflection $(x,y)\mapsto(x,-y)$, its trajectories are trajectories of the Hopf--Takens normal form studied in \cite{ZZ2005}, and they are Minkowski measurable \cite{ZZ2005}; see also \cite[Section~2]{HVZZ2023} and \cite[Theorem~2.1]{HRZ2025}. For $(m,n)\neq(1,1)$, the available results, namely \cite[Theorem~1]{HVZZ2023} for $m=n$ and Theorem~\ref{thm:main}, give Minkowski nondegeneracy only, and to the author's knowledge the Minkowski measurability of $\Gamma$ is open, already for $m=n\ge3$.
\end{remark}

\section{A scalar normal form and the decay of trajectories}\label{sec:nf}

This section is independent of Theorem~\ref{thm:main}. Its main result, Theorem~\ref{thm:decay}, gives the decay rate of $r(\varphi)$ along a trajectory; a box-dimension formula follows from it only under the additional hypothesis of Corollary~\ref{cor:dim}.

For a planar vector field
\[
X=(P,Q),
\]
set
\begin{equation}
\cT(X)=nxQ-myP,
\qquad
\cR(X)=2mx^{2m-1}P+2ny^{2n-1}Q,
\qquad
V=x^{2m}+y^{2n}.
\end{equation}

\begin{lemma}[Weighted polar identities]
In the weighted polar coordinates \eqref{eq:polar}, one has
\begin{equation}\label{eq:T}
\cT(X)=mn\,r^{m+n}\dot\varphi,
\end{equation}
and
\begin{equation}\label{eq:R}
\cR(X)=2mn\,r^{2mn-1}\dot r.
\end{equation}
Consequently,
\begin{equation}\label{eq:drdphi}
\frac{dr}{d\varphi}
=
\frac{r^{m+n+1}}{2V}\,\frac{\cR(X)}{\cT(X)}
\end{equation}
whenever $\cT(X)\neq 0$.
\end{lemma}

\begin{proof}
Write
\[
x=r^n\Cs(\varphi),
\qquad
y=r^m\Sn(\varphi).
\]
Differentiating and using \eqref{eq:cssn} gives
\begin{equation}
\begin{aligned}
P=\dot x&=n r^{n-1}\dot r\,\Cs(\varphi)-n r^n\Sn(\varphi)^{2n-1}\dot\varphi,\\
Q=\dot y&=m r^{m-1}\dot r\,\Sn(\varphi)+m r^m\Cs(\varphi)^{2m-1}\dot\varphi.
\end{aligned}
\end{equation}
For the angular part,
\begin{align*}
nxQ-myP
&=
n r^n\Cs(\varphi)\Bigl(m r^{m-1}\dot r\,\Sn(\varphi)+m r^m\Cs(\varphi)^{2m-1}\dot\varphi\Bigr)\\
&\quad
-m r^m\Sn(\varphi)\Bigl(n r^{n-1}\dot r\,\Cs(\varphi)-n r^n\Sn(\varphi)^{2n-1}\dot\varphi\Bigr)\\
&=
mn\,r^{m+n}\bigl(\Cs(\varphi)^{2m}+\Sn(\varphi)^{2n}\bigr)\dot\varphi\\
&=
mn\,r^{m+n}\dot\varphi,
\end{align*}
which proves \eqref{eq:T}.

For the radial part,
\begin{align*}
2mx^{2m-1}P
&=
2m r^{n(2m-1)}\Cs(\varphi)^{2m-1}
\Bigl(n r^{n-1}\dot r\,\Cs(\varphi)-n r^n\Sn(\varphi)^{2n-1}\dot\varphi\Bigr)\\
&=
2mn\,r^{2mn-1}\Cs(\varphi)^{2m}\dot r
-2mn\,r^{2mn}\Cs(\varphi)^{2m-1}\Sn(\varphi)^{2n-1}\dot\varphi,
\end{align*}
and similarly
\begin{align*}
2ny^{2n-1}Q
&=
2n r^{m(2n-1)}\Sn(\varphi)^{2n-1}
\Bigl(m r^{m-1}\dot r\,\Sn(\varphi)+m r^m\Cs(\varphi)^{2m-1}\dot\varphi\Bigr)\\
&=
2mn\,r^{2mn-1}\Sn(\varphi)^{2n}\dot r
+2mn\,r^{2mn}\Cs(\varphi)^{2m-1}\Sn(\varphi)^{2n-1}\dot\varphi.
\end{align*}
Adding the two identities, the mixed $\dot\varphi$ terms cancel and \eqref{eq:oval} gives
\[
2mx^{2m-1}P+2ny^{2n-1}Q
=
2mn\,r^{2mn-1}\bigl(\Cs(\varphi)^{2m}+\Sn(\varphi)^{2n}\bigr)\dot r
=
2mn\,r^{2mn-1}\dot r,
\]
which is \eqref{eq:R}. Finally, \eqref{eq:polar} implies
\[
V=x^{2m}+y^{2n}=r^{2mn},
\]
so dividing \eqref{eq:R} by \eqref{eq:T} yields
\[
\frac{\cR(X)}{\cT(X)}
=
\frac{2mn\,r^{2mn-1}\dot r}{mn\,r^{m+n}\dot\varphi}
=
2r^{2mn-m-n-1}\frac{dr}{d\varphi}
=
\frac{2V}{r^{m+n+1}}\frac{dr}{d\varphi},
\]
which is equivalent to \eqref{eq:drdphi}.
\end{proof}

By \eqref{eq:R}, $\cR(X)=\frac{d}{dt}V$ along trajectories, so it is natural to use
\[
z:=r^{2mn}=V
\]
as the radial variable; by \eqref{eq:T}, $dz/d\varphi=mn\,r^{m+n}\,\cR(X)/\cT(X)$ whenever $\cT(X)\neq0$. If along a small spiral trajectory one has a finite expansion
\begin{equation}\label{eq:rexp}
\frac{dr}{d\varphi}
=
r\sum_{j=1}^{N} h_j(\varphi)r^{2mnj}+O\bigl(r^{1+2mn(N+1)}\bigr)
\end{equation}
with continuous $T$-periodic coefficients $h_j$, then differentiation of $z=r^{2mn}$ yields
\begin{equation}
\frac{dz}{d\varphi}
=
\sum_{j=1}^{N} a_j(\varphi)z^{j+1}+O\bigl(z^{N+2}\bigr),
\qquad a_j(\varphi):=2mn\,h_j(\varphi).
\end{equation}
This is a scalar $T$-periodic equation, which can be simplified by the classical method of periodic averaging with near-identity transformations, the size of $z$ playing the role of the small parameter; see \cite[Chapters~2 and~3]{Sanders2007}, where the nonuniqueness of the averaging transformation is also discussed. For definiteness we fix the scheme in which all periodic primitives have zero mean, and we give the argument in full. For a $T$-periodic function $f$ we write $\mean{f}:=\frac1T\int_0^Tf(\tau)\,d\tau$. Throughout this section, $O(r^{1+2mn(N+1)})$, $O(z^{N+2})$ and $O(\rho^{N+2})$ denote remainders bounded by $C|r|^{1+2mn(N+1)}$, $C|z|^{N+2}$ and $C|\rho|^{N+2}$, respectively, uniformly in $\varphi$, for sufficiently small $|r|$, $|z|$, $|\rho|$.

\begin{lemma}[A uniform domain for near-identity substitutions]\label{lem:domain}
Fix $N\ge 1$. For each $\nu\in\{0,\dots,N-1\}$ let $v_{\nu+1}$ be a bounded $T$-periodic function and set
\[
\Phi_\nu(\varphi,\xi):=\xi+v_{\nu+1}(\varphi)\xi^{\nu+2}.
\]
Then there exists $\delta_*>0$ such that the following hold.
\begin{enumerate}[label=\textnormal{(\roman*)},leftmargin=2em]
\item For every $\nu\in\{0,\dots,N-1\}$ and every fixed $\varphi$, the map $\xi\mapsto \Phi_\nu(\varphi,\xi)$ is a $C^1$ diffeomorphism from $(-\delta_*,\delta_*)$ onto its image.
\item There exist constants $0<c_\nu\le C_\nu$ such that for all $|\xi|<\delta_*$ and all $\varphi$,
\[
c_\nu |\xi|\le |\Phi_\nu(\varphi,\xi)|\le C_\nu |\xi|.
\]
\item If $|\rho|<\delta$ with
\[
\delta:=\frac{\delta_*}{\prod_{\nu=0}^{N-1} C_\nu},
\]
then every intermediate output in the finite composition
\[
\Phi_0\bigl(\varphi,\Phi_1(\varphi,\cdots \Phi_{N-1}(\varphi,\rho)\cdots)\bigr)
\]
remains in $(-\delta_*,\delta_*)$, uniformly in $\varphi$.
\end{enumerate}
\end{lemma}

\begin{proof}
Since only finitely many bounded functions are involved,
\[
M:=\max_{0\le \nu\le N-1}\ \sup_{\varphi\in\mathbb R}|v_{\nu+1}(\varphi)|<\infty.
\]
Choose $\delta_*>0$ so small that
\[
(\nu+2)M\delta_*^{\nu+1}\le \frac12
\qquad\text{for every }\nu\in\{0,\dots,N-1\}.
\]
Then, for $|\xi|<\delta_*$,
\[
\partial_\xi\Phi_\nu(\varphi,\xi)=1+(\nu+2)v_{\nu+1}(\varphi)\xi^{\nu+1}
\]
satisfies
\[
\frac12\le \partial_\xi\Phi_\nu(\varphi,\xi)\le\frac32
\]
uniformly in $\varphi$. Hence each map is strictly increasing and is a $C^1$ diffeomorphism from $(-\delta_*,\delta_*)$ onto its image. Moreover, since $M\delta_*^{\nu+1}\le 1/(2(\nu+2))\le 1/4$,
\[
|\Phi_\nu(\varphi,\xi)-\xi|
\le M|\xi|^{\nu+2}
\le M\delta_*^{\nu+1}|\xi|
\le \frac14|\xi|,
\]
so (ii) holds, for example, with $c_\nu=1/2$ and $C_\nu=3/2$.

For the finite composition, set
\[
\rho_N:=\rho,
\qquad
\rho_\nu:=\Phi_\nu(\varphi,\rho_{\nu+1})
\quad (\nu=N-1,N-2,\dots,0).
\]
Applying the upper estimate in (ii) repeatedly, and using $C_j\ge1$, we obtain
\[
|\rho_\nu|
\le\Bigl(\prod_{j=\nu}^{N-1}C_j\Bigr)|\rho|
\le\Bigl(\prod_{j=0}^{N-1}C_j\Bigr)|\rho|
<\delta_*
\]
for every $\nu\in\{0,\dots,N\}$ whenever $|\rho|<\delta$. In particular, every input to the next map lies in $(-\delta_*,\delta_*)$, uniformly in $\varphi$. This proves (iii).
\end{proof}

\begin{theorem}[Finite-order normal form]\label{thm:nf}
Fix $N\ge 1$. Consider a scalar equation
\begin{equation}\label{eq:zeq}
\frac{dz}{d\varphi}
=
\sum_{j=1}^{N} a_j(\varphi)z^{j+1}+O\bigl(z^{N+2}\bigr),
\end{equation}
where each $a_j$ is continuous and $T$-periodic. Then there exist constants $g_1,\dots,g_N$ and a change of variables
\begin{equation}
z=\Psi(\varphi,\rho),
\end{equation}
defined for $|\rho|$ sufficiently small, with the following properties:
\begin{enumerate}[label=\textnormal{(\roman*)}]
\item for each fixed $\varphi$, the map $\rho\mapsto \Psi(\varphi,\rho)$ is a $C^1$ diffeomorphism from a neighbourhood of $0$ onto its image;
\item $\Psi$ is $T$-periodic in $\varphi$ and has the expansion
\begin{equation}
\Psi(\varphi,\rho)
=
\rho+u_1(\varphi)\rho^2+\cdots+u_N(\varphi)\rho^{N+1}+O\bigl(\rho^{N+2}\bigr),
\end{equation}
where $u_1,\dots,u_N$ are $C^1$ and $T$-periodic;
\item under the change $z=\Psi(\varphi,\rho)$, equation \eqref{eq:zeq} is transformed into
\begin{equation}\label{eq:rhoeq}
\frac{d\rho}{d\varphi}
=
\sum_{j=1}^{N} g_j\rho^{j+1}+O\bigl(\rho^{N+2}\bigr).
\end{equation}
\end{enumerate}
Here (iii) means the following: if $z$ is a $C^1$ function with $|z(\varphi)|$ sufficiently small that satisfies \eqref{eq:zeq}, then the $C^1$ function $\rho$ defined by $z(\varphi)=\Psi(\varphi,\rho(\varphi))$ satisfies \eqref{eq:rhoeq}. Moreover, there exists $\delta>0$ such that for every $\varphi$ the map $\rho\mapsto \Psi(\varphi,\rho)$ is a $C^1$ diffeomorphism from $(-\delta,\delta)$ onto its image. In the construction, the periodic primitive $v_{\nu+1}$ at each step is normalized by $\mean{v_{\nu+1}}=0$; the constants $g_j$ are those produced by this zero-mean scheme.
\end{theorem}

\begin{proof}
We construct the normal form inductively. For $\nu\in\{0,1,\dots,N\}$ we seek a variable $z_\nu$ satisfying
\begin{equation}\label{eq:stepnu}
\frac{dz_\nu}{d\varphi}
=
\sum_{j=1}^{\nu}g_jz_\nu^{j+1}
+
\sum_{j=\nu+1}^{N}a_j^{(\nu)}(\varphi)z_\nu^{j+1}
+
O(z_\nu^{N+2}),
\end{equation}
where the coefficient functions are continuous and $T$-periodic. For $\nu=0$ this is the given equation with $z_0=z$ and $a_j^{(0)}=a_j$.

Assume that \eqref{eq:stepnu} holds for some $\nu<N$. Define
\[
g_{\nu+1}:=\mean{a_{\nu+1}^{(\nu)}}
=\frac1T\int_0^T a_{\nu+1}^{(\nu)}(\tau)\,d\tau.
\]
The function $a_{\nu+1}^{(\nu)}-g_{\nu+1}$ has zero mean, so the equation
\[
v_{\nu+1}'(\varphi)=a_{\nu+1}^{(\nu)}(\varphi)-g_{\nu+1}
\]
has a $T$-periodic $C^1$ solution, which becomes unique if we impose $\mean{v_{\nu+1}}=0$. We use the substitution
\begin{equation}\label{eq:subst}
z_\nu=z_{\nu+1}+v_{\nu+1}(\varphi)z_{\nu+1}^{\nu+2}.
\end{equation}
By Lemma~\ref{lem:domain}(i),(ii), applied to this single map, \eqref{eq:subst} is a $C^1$ diffeomorphism for sufficiently small $|z_{\nu+1}|$, uniformly in $\varphi$, and $z_\nu\asymp z_{\nu+1}$.

Write $w=z_{\nu+1}$, $v=v_{\nu+1}$, and
\[
B(\varphi,w):=(\nu+2)v(\varphi)w^{\nu+1}.
\]
Along a solution, $w$ is $C^1$ by the implicit function theorem, because $\Phi_\nu$ is $C^1$ in $(\varphi,\xi)$ and $\partial_\xi\Phi_\nu\ge1/2$. Differentiating \eqref{eq:subst} gives
\[
(1+B(\varphi,w))w'
=
F_\nu\bigl(\varphi,w+v(\varphi)w^{\nu+2}\bigr)
-v'(\varphi)w^{\nu+2},
\]
where $F_\nu$ denotes the right-hand side of \eqref{eq:stepnu}. On the neighbourhood given by Lemma~\ref{lem:domain}, $|B|\le1/2$, and hence
\begin{equation}\label{eq:quot}
w'
=
\frac{F_\nu(\varphi,w+v(\varphi)w^{\nu+2})-v'(\varphi)w^{\nu+2}}
{1+B(\varphi,w)}.
\end{equation}
We expand the right-hand side up to degree $N+1$ in $w$.

Substitution into the polynomial part of $F_\nu$ and the binomial formula give continuous $T$-periodic functions $A_j(\varphi)$ such that the numerator in \eqref{eq:quot} equals
\begin{align*}
&\sum_{j=1}^{\nu}g_jw^{j+1}
+\bigl(a_{\nu+1}^{(\nu)}(\varphi)-v'(\varphi)\bigr)w^{\nu+2}
+\sum_{j=\nu+2}^{N}A_j(\varphi)w^{j+1}
+O(w^{N+2}).
\end{align*}
The remainder of $F_\nu$ is still $O(w^{N+2})$ because $z_\nu\asymp w$. To divide by $1+B$ without losing higher-order terms, we use the identity
\[
\frac1{1+B}=\sum_{i=0}^{N}(-B)^i+\frac{(-B)^{N+1}}{1+B}.
\]
Since $B=O(w^{\nu+1})$ and the numerator is $O(w^2)$, the last term contributes
\[
O(w^2)\,O(w^{(N+1)(\nu+1)})\subset O(w^{N+3}),
\]
which is contained in $O(w^{N+2})$. Multiplying the finite sum by the numerator and collecting powers up to degree $N+1$, we obtain continuous $T$-periodic coefficients $a_j^{(\nu+1)}$ such that
\[
w'
=
\sum_{j=1}^{\nu}g_jw^{j+1}
+\bigl(a_{\nu+1}^{(\nu)}-v'\bigr)w^{\nu+2}
+\sum_{j=\nu+2}^{N}a_j^{(\nu+1)}(\varphi)w^{j+1}
+O(w^{N+2}).
\]
The factor $1/(1+B)$ does not change the terms of degree at most $\nu+2$: its first nonconstant term has degree $\nu+1$, while the numerator starts at degree $2$, so their first product has degree $\nu+3$. Hence the coefficients $g_1,\dots,g_\nu$ are unchanged. By the definition of $v$,
\[
a_{\nu+1}^{(\nu)}-v'=g_{\nu+1},
\]
which completes the induction step.

After $N$ steps, set $\rho=z_N$. The total change of variables is the finite composition
\[
\Psi(\varphi,\rho)
=
\Phi_0\!\bigl(\varphi,\Phi_1(\varphi,\ldots,
\Phi_{N-1}(\varphi,\rho)\ldots)\bigr),
\qquad
\Phi_\nu(\varphi,\xi)=\xi+v_{\nu+1}(\varphi)\xi^{\nu+2}.
\]
Lemma~\ref{lem:domain} gives a common domain on which all intermediate values are defined. Substituting and collecting powers gives
\[
\Psi(\varphi,\rho)
=
\rho+u_1(\varphi)\rho^2+\cdots+u_N(\varphi)\rho^{N+1}
+O(\rho^{N+2}),
\]
with $C^1$, $T$-periodic coefficient functions, and the transformed equation is
\[
\rho'=g_1\rho^2+\cdots+g_N\rho^{N+1}+O(\rho^{N+2}).
\]
Finally, $\partial_\rho\Psi(\varphi,0)=1$ for every $\varphi$. By continuity of this derivative and compactness of $[0,T]$, there is $\delta>0$ such that
\[
|\partial_\rho\Psi(\varphi,\rho)-1|\le\frac12
\qquad (\varphi\in[0,T],\ |\rho|\le\delta).
\]
Thus $\partial_\rho\Psi\ge1/2$, and by periodicity each map $\rho\mapsto\Psi(\varphi,\rho)$, $\varphi\in\R$, is a $C^1$ diffeomorphism from $(-\delta,\delta)$ onto its image.
\end{proof}

\begin{definition}[Scalar normal-form coefficients]
The constants $g_1,\dots,g_N$ of Theorem~\ref{thm:nf} are called the \emph{scalar normal-form coefficients} of \eqref{eq:zeq} up to order $N$ (with respect to the zero-mean scheme).
\end{definition}

\begin{remark}[Normal-form coefficients and averages]\label{rem:averages}
In general the coefficients $g_j$ are not the averages $\mean{a_j}$, and hence not $2mn$ times the averages of the coefficients $h_j$ in \eqref{eq:rexp}. For $j\le2$ they are. Indeed, $g_1=\mean{a_1}$ by construction. Starting from
\[
z'=a_1(\varphi)z^2+a_2(\varphi)z^3+O(z^4),
\]
the first step uses
\[
z=\rho+v_1(\varphi)\rho^2,
\qquad
v_1'(\varphi)=a_1(\varphi)-g_1,
\qquad
g_1=\mean{a_1},
\]
and collecting cubic terms gives
\[
\rho'
=
g_1\rho^2+
\bigl(a_2(\varphi)+2v_1(\varphi)(a_1(\varphi)-g_1)\bigr)\rho^3
+
O(\rho^4).
\]
Since $2v_1(a_1-g_1)=2v_1v_1'=(v_1^2)'$ is the derivative of a $T$-periodic function,
\[
g_2=\mean{a_2+(v_1^2)'}=\mean{a_2}.
\]
From order $3$ on, lower-order coefficients contribute. If $z'=a_1z^2+a_2z^3+a_3z^4+O(z^5)$, the first step produces the quartic coefficient $a_3+a_2v_1+g_1v_1^2-(v_1^3)'$, the second step (with $v_2$ of zero mean) subtracts $g_1v_2$ from it, and therefore
\[
g_3=\mean{a_3}+\mean{a_2v_1}+g_1\mean{v_1^2}.
\]
For instance, for $T=2\pi$, $a_1=\cos\varphi$, $a_2=-\sin\varphi$ and $a_3=0$, all averages $\mean{a_j}$ vanish, whereas $v_1=\sin\varphi$, $g_1=g_2=0$ and $g_3=-\mean{\sin^2\varphi}=-1/2<0$, so Theorem~\ref{thm:decay} below applies with $\ell=3$. For this reason the hypotheses of Theorem~\ref{thm:decay} are formulated in terms of the coefficients $g_j$ and not in terms of averages.
\end{remark}

\begin{theorem}[Decay rate]\label{thm:decay}
Assume that along a small spiral trajectory the function $z=r^{2mn}$ satisfies an equation of the form \eqref{eq:zeq} with continuous $T$-periodic coefficients $a_1,\dots,a_N$, that is, $dz/d\varphi-\sum_{j=1}^N a_j(\varphi)z^{j+1}=O(z^{N+2})$ along the trajectory. Let $g_1,\dots,g_N$ be its scalar normal-form coefficients, and let $z=\Psi(\varphi,\rho)$ be the change of variables of Theorem~\ref{thm:nf}; for $|\rho|$ sufficiently small, $\Psi$ preserves sign, so that positive $z$ corresponds to positive $\rho$. Suppose moreover that for some $\ell\in\{1,\dots,N\}$,
\begin{equation}\label{eq:gl}
g_1=\cdots=g_{\ell-1}=0,
\qquad
g_{\ell}<0.
\end{equation}
Then, for all sufficiently large $\varphi$, the function $\rho(\varphi)$ defined by $z(\varphi)=\Psi(\varphi,\rho(\varphi))$ with $\rho(\varphi)$ small is positive, and
\begin{equation}
z(\varphi)\asymp \varphi^{-1/\ell},
\qquad
r(\varphi)\asymp \varphi^{-1/(2mn\ell)}.
\end{equation}
\end{theorem}

\begin{proof}
By Theorem~\ref{thm:nf} and \eqref{eq:gl}, the change $z=\Psi(\varphi,\rho)$ transforms \eqref{eq:zeq} into
\[
\frac{d\rho}{d\varphi}
=
g_{\ell}\rho^{\ell+1}+O\bigl(\rho^{\ell+2}\bigr),
\qquad g_{\ell}<0.
\]
Moreover, the proof of Theorem~\ref{thm:nf} gives $\delta>0$ such that for every $\varphi$ the map $\rho\mapsto \Psi(\varphi,\rho)$ is a $C^1$ diffeomorphism on $(-\delta,\delta)$ and
\[
|\partial_\rho\Psi(\varphi,\rho)-1|\le \frac12
\qquad\text{for all }\varphi\in\R,\ |\rho|\le \delta.
\]
Hence, for each fixed $\varphi$, the map $\rho\mapsto \Psi(\varphi,\rho)$ is strictly increasing on $(-\delta,\delta)$ and satisfies $\Psi(\varphi,0)=0$, so it preserves sign there. By the mean value theorem and the derivative bound, the image $J_\varphi:=\Psi(\varphi,(-\delta,\delta))$ contains $(-\delta/2,\delta/2)$ for every $\varphi$. Since $|\partial_\rho\Psi-1|\le1/2$, the mean value theorem also gives $|\Psi(\varphi,\rho)-\rho|\le|\rho|/2$, that is, for $z=\Psi(\varphi,\rho)$,
\[
\frac12|\rho|\le |z|\le \frac32|\rho|
\qquad\text{for all }|\rho|<\delta\text{ and all }\varphi.
\]

Along the given small spiral trajectory, $z=r^{2mn}>0$ and $z(\varphi)\to0$ as $\varphi\to\infty$. Choose $\varphi_1\ge1$ so large that $0<z(\varphi)<\delta/2$ for all $\varphi\ge \varphi_1$; then $z(\varphi)\in J_\varphi$, and we may define
\[
\rho(\varphi):=\Psi(\varphi,\cdot)^{-1}\bigl(z(\varphi)\bigr)\in(-\delta,\delta)
\qquad(\varphi\ge\varphi_1).
\]
Since the functions $v_{\nu+1}$ in the construction of $\Psi$ are $C^1$, $\Psi$ is $C^1$ in $(\varphi,\rho)$; as $\partial_\rho\Psi\ge1/2$, the implicit function theorem shows that $\rho$ is $C^1$. Since $z(\varphi)=\Psi\bigl(\varphi,\rho(\varphi)\bigr)$, the function $\rho$ satisfies the normalized equation
\[
\frac{d\rho}{d\varphi}=g_{\ell}\rho^{\ell+1}+O\bigl(\rho^{\ell+2}\bigr)
\qquad (\varphi\ge \varphi_1).
\]
Since $\Psi(\varphi,\cdot)$ preserves sign, $\rho(\varphi)>0$, and $\frac23 z(\varphi)\le\rho(\varphi)\le 2z(\varphi)$; in particular $\rho(\varphi)\to0$.

Choose $\eta$ with $0<\eta<|g_{\ell}|/2$. For sufficiently small positive $\rho$,
\[
-(|g_{\ell}|+\eta)\rho^{\ell+1}
\le
\frac{d\rho}{d\varphi}
\le
-(|g_{\ell}|-\eta)\rho^{\ell+1},
\]
and after increasing $\varphi_1$ we may assume that these inequalities hold along the trajectory for all $\varphi\ge \varphi_1$. Since $\rho(\varphi)>0$,
\[
\frac{d}{d\varphi}(\rho^{-\ell})=-\ell\rho^{-\ell-1}\frac{d\rho}{d\varphi},
\]
and therefore
\[
\ell(|g_{\ell}|-\eta)
\le
\frac{d}{d\varphi}(\rho^{-\ell})
\le
\ell(|g_{\ell}|+\eta)
\qquad(\varphi\ge\varphi_1).
\]
Integrating from $\varphi_1$ to $\varphi$ gives
\[
\rho(\varphi_1)^{-\ell}
+\ell(|g_{\ell}|-\eta)(\varphi-\varphi_1)
\le
\rho(\varphi)^{-\ell}
\le
\rho(\varphi_1)^{-\ell}
+\ell(|g_{\ell}|+\eta)(\varphi-\varphi_1).
\]
For $\varphi\ge2\varphi_1$ the left-hand side is at least $\tfrac12\ell(|g_\ell|-\eta)\varphi$, and, since $\varphi\ge1$, the right-hand side is at most $\bigl(\rho(\varphi_1)^{-\ell}+\ell(|g_\ell|+\eta)\bigr)\varphi$. Hence $\rho(\varphi)^{-\ell}\asymp\varphi$, that is,
\[
\rho(\varphi)\asymp \varphi^{-1/\ell}.
\]
Since $z(\varphi)\asymp \rho(\varphi)$, it follows that $z(\varphi)\asymp \varphi^{-1/\ell}$, and since $z=r^{2mn}$, we obtain $r(\varphi)\asymp \varphi^{-1/(2mn\ell)}$.
\end{proof}

\begin{corollary}[Box dimension under a comparison hypothesis]\label{cor:dim}
Assume the hypotheses of Theorem~\ref{thm:decay}. If, in addition, $m\ge n$ and the spiral tail satisfies the hypothesis of Definition~\ref{def:comparison} with
\[
p=\frac{1}{2m\ell},
\qquad
q=\frac{1}{2n\ell},
\]
then
\[
\dimB\Gamma=2-\frac{1+\frac{n}{m}}{1+2n\ell}.
\]
\end{corollary}

\begin{proof}
By Theorem~\ref{thm:decay}, $r(\varphi)\asymp\varphi^{-1/(2mn\ell)}$, so $|x(\varphi)|\le r^n\le C\varphi^{-p}$ and $|y(\varphi)|\le r^m\le C\varphi^{-q}$; this is why $S_{p,q}$ with these exponents is the natural comparison spiral. Since $\ell\ge 1$, we have $p\le 1/2<1$, and since $m\ge n$, we have $p\le q$. Hence Lemma~\ref{lem:transfer} applies and gives
\[
\dimB\Gamma=2-\frac{p+q}{1+q}=2-\frac{1+\frac{n}{m}}{1+2n\ell}.
\]
The conclusion thus follows from the comparison hypothesis and Lemma~\ref{lem:transfer}; Theorem~\ref{thm:decay} only identifies the exponents $p$ and $q$.
\end{proof}

\section{Concluding remarks}\label{sec:concl}

Theorem~\ref{thm:main} and Corollary~\ref{cor:unstable} prove the conjecture of \cite{HVZZ2023} for all odd $m\ge n\ge1$ and $k\ge1$: every spiral trajectory of \eqref{eq:system}, and of the field with the opposite signs, is Minkowski nondegenerate and has box dimension $2-(1+n/m)/(1+2nk)$; by Remark~\ref{rem:swap}, the case $m<n$ follows by exchanging the roles of $m$ and $n$. Following the scheme of \cite{BFF2022}, the proof uses the explicit form of the trajectories, in particular the relation \eqref{eq:turns} between consecutive turns: the lengths of the turns give the upper bound, and the separation of consecutive turns in a fixed sector of the first quadrant gives a matching lower bound. The author does not know whether the trajectories of \eqref{eq:system} are bi-Lipschitz equivalent to the elliptical polynomial spirals $S_{p,q}$.

For a general vector field with a degenerate focus of type $(m,n)$, Theorem~\ref{thm:decay} gives the decay rate of $r(\varphi)$ once a scalar equation of the form \eqref{eq:zeq} is known along a trajectory. To deduce the box dimension from it, we need additional geometric information, since a relation such as \eqref{eq:turns} is not part of its hypotheses; Corollary~\ref{cor:dim} obtains the box dimension under a bi-Lipschitz comparison hypothesis. This leaves two natural problems: to find conditions on the vector field which guarantee \eqref{eq:zeq} along its trajectories, and to find geometric conditions under which trajectories without an explicit turn structure satisfy the conclusion of Theorem~\ref{thm:main}, either through a transfer lemma such as Lemma~\ref{lem:transfer} or through a direct covering argument. The Minkowski measurability of the trajectories of \eqref{eq:system} for $(m,n)\neq(1,1)$ also remains open (Remark~\ref{rem:measurability}).

\section*{Acknowledgments}
Between 1~April and 25~November 2026, the author used ChatGPT (OpenAI; GPT-5.4 Pro, GPT-5.5 Pro, GPT-5.6 Sol Pro, GPT-6 Pro/GPT-6 Astra, and GPT-6.2 Astra) and Claude (Anthropic; Opus~4.6, 4.7, 4.8, 5, and 5.5; Fable~5, 5.1, and 5.2), adopting successive versions as they became available. These tools assisted with proof checking and revision, language editing, and \LaTeX\ preparation. The author reviewed and edited the AI-assisted material and accepts full responsibility for the manuscript.

\section*{Use of generative artificial intelligence}
The AI tools identified in the Acknowledgments assisted with examining mathematical arguments for correctness and logical consistency, revising proofs and their exposition, improving the English text, and editing the \LaTeX\ source. The author assessed the suggestions and made the final decisions concerning the mathematical content and its presentation.

\end{document}